\documentclass[10pt,a4paper,oneside]{article} 
\usepackage{amssymb,amsmath} 
\usepackage[all]{xy}
\usepackage{fancyheadings} 
\usepackage{latexsym}
\usepackage{amsthm}

\newtheorem{theo}{{\bf Theorem}}[subsection]
\newtheorem{proposition}{{\bf Proposition}}[subsection]
\newtheorem{corollaire}{{\bf Corollary}}[subsection]
\newtheorem{lemme}{{\bf Lemma}}[subsection]

\newtheorem{remarques}{{\bf Remarks}}[subsection]

\newcommand{\Cal}{\mathcal} 
\newtheorem{definition}{{\bf Definition}}[subsection]

\begin{document}

{ \bf  CLASSICAL YANG-BAXTER EQUATION AND\\ LEFT INVARIANT AFFINE GEOMETRY ON LIE GROUPS}

\vskip 1truecm
 ANDRE DIATTA\footnote{ \footnotesize The first author was partially supported by Enterprise Ireland.} 
 AND ALBERTO MEDINA
\vskip 1truecm


\begin{abstract}{ \footnotesize

Let $G$ be a Lie group with  Lie algebra $ \Cal G: = T_\epsilon G$ and $T^*G = \Cal G^* \rtimes G$ its cotangent bundle considered as a Lie group, where $G$ acts on $\Cal G^*$ via the coadjoint action.
 We show that there is a 1-1 correspondance between the skew-symmetric solutions $r\in \wedge^2 \Cal G$ of the Classical Yang-Baxter Equation in $G$, and the set of connected Lie subgroups of  $T^*G$ which carry a left invariant affine structure and whose Lie algebras are  lagrangian graphs in $ \Cal G \oplus \Cal G^*$.

An invertible solution $r$ endows $G$ with a left invariant symplectic structure and hence a left invariant affine structure. In this case we prove that the Poisson Lie tensor $\pi := r^+ - r^-$ is polynomial of degree at most 2 and the double Lie groups of $(G,\pi)$  also carry a canonical left invariant affine structure.

In the general case of (non necessarly invertible) solutions $r$, we supply a necessary and suffisant condition to the geodesic completness of the associated affine structure.}
 
\end{abstract}

\section{\bf Introduction-Summary}
Let  $G$ be a Lie group whose Lie algebra is denoted by $\Cal G:= T_\epsilon G$, where $\epsilon$ stands for the unit of $G$.

Since the end of $60$s-the early $70$s (see e.g. {\bf \cite{Chu}}), it's well known that if $w^+$ is a left invariant symplectic form on $G$, the formula 
 \begin{equation}\label{PSG}
 \omega^+ (\nabla_{x^+}y^+, z^+):= -\omega^+ (y^+, [x^+,z^+])
\end{equation}
defines a (locally) flat and torsion free connection $\nabla$ (i.e an affine structure) in $G$ which is left invariant. Here  $x^+$ is the left invariant vector field satisfying $x^+_\epsilon = x$.
Such a formula plays a crucial role in the study of symplectic and k\"ahlerian Lie groups developed in  {\bf \cite{Li-Me}}, {\bf \cite{Li-Me88}}, {\bf \cite{Lichne91}}, {\bf \cite{Me-Re91}}, {\bf \cite{Da-Me1}}, {\bf \cite{Da-Me2}}.

The problem of finding those Lie groups (necessary solvable) which admit a
geodesically complete, left invariant flat torsion free connection, is an open
problem (see Milnor {\bf \cite{Milnor}}), and few such Lie groups are known.

 Left invariant affine structures on $G$ bijectively correspond to left symmetric algebra (LSA, also called Koszul-Vinberg algebra) structures on the Lie
algebra  $\Cal G$ of $G$, compatible with bracket of  $\Cal G$ ( J.L. Koszul
{\bf \cite{Kosz55}}, {\bf \cite{Kos61}}).
  Such structures first arised in the works of J.L. Koszul {\bf \cite{Kosz55}}, {\bf \cite{Kos61}}, E.B. Vinberg {\bf \cite{Vinb}}, Y. Matsushima.
 Left symmetric algebras has been studied in Helmstetter, Medina, Dardi\'e, Chu in {\bf \cite{Hel79}}, {\bf\cite{Me81}}, {\bf \cite{Da-Me1}}, {\bf \cite{Da-Me2}},  {\bf \cite{Chu}}
... etc.

Later, M. Bordemann ({\bf \cite{Bor90}}) shows that a solution
of the Classical Yang Baxter Equation (CYBE) $r\in \wedge^2 \Cal G$  always
determines a left invariant affine structure on any Lie group $G^*$ whose Lie
algebra is the dual $\Cal G^*$  (relative to $r$) of $\Cal G$. The corresponding left invariant connection on $ G^*$, is given by the formula
\begin{equation}\label{PSGDual0}
  \nabla_{\alpha^+}\beta^+:=
(ad_{r(\alpha)}^*\beta)^+.
 \end{equation} 
Let's identify the cotangent
bundle $T^*G$ of $G$ and the trivial bundle $\Cal G^*\times G$, by means of
left translations. And let's endow this latter with a structure of
semi-direct product group $\Cal G^*\rtimes G$ of $G$ and the abelian Lie group
$\Cal G^*$, where $G$ acts in $\Cal G^*$ by the coadjoint action.
 
 Denote
$\pi:= r^+-r^-$ the Poisson Lie tensor associated to $r$, and  $\Cal D(\Cal
G)$ the double Lie algebra of $(G,\pi)$. Here are some of the results we prove
in this work: 
\\
\\
 - there is a bijective correspondance between the
skew-symmetric solutions of the Classical Yang Baxter Equation (CYBE) on $\Cal
G := Lie(G)$  and the connected Lie subgroups of $T^*G$ equipped with left
invariant affine structures whose Lie algebras are lagrangian graphs of a
$\Cal G $-valued skew-symmetric linear mapping of $\Cal G^*$. Here  the Lie
algebra $t^*\Cal G:= \Cal G^*\rtimes\Cal G$ of $T^*G$, is endowed with its
canonical hyperbolic orthogonal structure ie its canonical dual pairing (theorem {\bf \ref{dualYBC-Lag} }).
\\
 - let  $r\in \wedge^2 \Cal G$ be a solution of the CYBE, then the affine
structure  defined by {\bf (\ref{PSGDual0})} in any dual Lie group $G^*$,
 is geodesically
complete if and only if the symplectic leaf through $\epsilon$, of the left
invariant Poisson structure $r^+$, is a unimodular (symplectic) Lie group
(theorem {\bf \ref{GeoComplAf}}). Thus, the geodesic completness of  $G^*$
implies the solvability of $G^*$, thanks to a result due to
Lichnerowicz-Medina {\bf \cite{Li-Me}}.
\\
\\
  In the case where $r$ is an
invertible solution of the CYBE, we prove: 
\\
  {\bf -} the existence of left
invariant affine structure, defined by $r$, on every connected Lie group $D(
G)$ whose Lie algebra is the double  Lie algebra $\Cal D (\Cal G)$ (theorem
{\bf \ref{doubleYBC3}}). Such a structure is geodesically complete if and only
if $G$  is unimodular (and solvable) ( see theorem {\bf\ref{GeoComplDoubl}}).; 
\\
 {\bf -} the existence on $ D (G)$ of a  left
invariant complex structure given by $r$ (proposition {\bf \ref{doubleYBC1}});
\\
 {\bf -} that the Poisson Lie tensor $\pi = r^+ - r^-$ is polynomial of
degree $2$ (theorem {\bf \ref{local-aff1}}).

  Moreover, if the differential
2-form $\omega^+$ is exact, we prove that the symplectic leaves of $\pi$ are
coverings over coadjoint orbits of $G$ (Corollary {\bf \ref{feuille-coadj} }). 
\\

We would like to thank prof. D. V. Alekseevsky for various hepful discussions and comments.
\\
\\
 {\bf Reminders and notations.}

 An affine structure on a manifold $M$, is
given by a maximal atlas of charts (affine atlas) to an affine space with
transition functions extending to affine transformations. Equivalently, an
affine stucture on $M$ is given by  an immersion, a so called developing map, 
from the universal cover $\tilde M$ of $M$ to an affine space (e.g $\Bbb R^n$
), equivariant with respect to the so called holonomy representation  $h:
\pi_1 (M)\to Aff(\Bbb R^n)$, where $\pi_1(M)$ is the group of deck
transformations of $\tilde M$.

 A morphism of affine manifolds is a map
whose expression  in affine atlas extends to an affine map. An affine stucture
in a Lie group $G$ is said to be left invariant if left translations are
morphisms of the affine manifold $G$.  

 A Poisson manifold $M$ is a
manifold whose ring of functions $\Cal C^\infty (M)$ carries a Lie algebra
bracket $\{,\}$ which is a first order linear differential operator, in each
of its arguments.
 The Hamiltonian vector field associated to an $f\in \Cal
C^\infty (M)$ is $X_f:= \{f,\hskip 1truemm .\}=:\Lambda (df, \hskip 1truemm . ) =: \Lambda^\sharp
(df)$,  where $\Lambda$ is the associated bivector field. The tensor $\Lambda$
must satisfy $[\Lambda,\Lambda] = 0$, here $[,]$ is the so-called Schouten
bracket ( see e.g {\bf \cite{Lichne77}}).

 Every Poisson manifold foliates
into symplectic manifolds (its symplectic leaves), which are the maximal
integrals of the involutive distribution spanned by all the Hamiltonian vector
fields. In general, it's a difficult problem to describe symplectic leaves of
a Poisson manifold. One of the most important families of Poisson manifolds is the one of Poisson Lie groups. 

 A Poisson Lie group $(G,\pi)$ is a Lie group $G$
together with a Poisson structure $\pi$ such that the multiplication is a
Poisson morphism. Set $[\alpha,\beta]_*:= d_\epsilon \{f,g\}$, where $\alpha:=
d_\epsilon f$ and $\beta:= d_\epsilon g\in \Cal G^*$  are the respective
differentials of $f,g\in \Cal C^\infty (G)$ at the unit $\epsilon$. Thus
$(\Cal G^*, [,]_*)$ is a Lie algebra called the dual of $(G,\pi)$ (or the dual
of $\Cal G$). While the double of  $(G,\pi)$ is the vector space $\Cal
G\oplus\Cal G^* $ equipped with the Lie bracket $[(x,\alpha), (y,\beta) ] :=(
[x,y]_{\Cal G} + ad_\alpha^* y - ad_\beta^*x,[\alpha,\beta]_* + ad^*_x\beta - 
ad^*_y\alpha) $. Here both coadjoint actions  of $\Cal G$ on  $\Cal G^*$ and
$(\Cal G^*,[,]_*)$ on $\Cal G$, are considered (see {\bf \cite{Dr83}}).

Denoting $r^+$ a left invariant bivector field in $G$, whose value at the
unit is $r^+_\epsilon=:r$ and $[r^+,r^+]$ the Schouten bracket of $r^+$ and itself, the equation
 \begin{equation}\label{CYBE}
 r\in\wedge^2 \Cal G, \hskip 2truecm [r,r]:= [r^+,r^+]_\epsilon =0
 \end{equation}
is called the Classical Yang Baxter Equation (CYBE) in $Lie(G)=: \Cal G$, (see
{\bf \cite{Dr83}}).  For a solution $r$ of the CYBE, the tensor $\pi:
=r^+-r^- $ defines a Poisson Lie structure in $G$ (see also {\bf \cite{Lu-We}
}, {\bf \cite{STS83} }, ...).

 We'll identify $\otimes^2\Cal G$ and the space of linear maps $\Cal G^*\to \Cal G$ and will denote an element $r$ of $\otimes^2\Cal G$ and the corresponding linear map $ \alpha\mapsto r^\sharp(\alpha):=r(\alpha,.) $ by the same symbol.  With such identifications, the 3-vector $[r,r]$ above is given by $[r,r] (\alpha,\beta,\gamma):= <[r(\alpha),r(\beta)], \gamma>+<[r(\beta),r(\gamma)],\alpha>+<[r(\gamma),r(\alpha)],\beta>$ for every $\alpha,\beta,\gamma\in \Cal G^*$, where $[,]$ is the Lie bracket in $\Cal G$ and $<X,\alpha>:= \alpha(X)$ is the usual canonical  dual pairing between vectors and linear forms on $\Cal G$.

\begin{definition}

 A Lie algebra $\Cal G$ will be termed a Manin algebra if it has a (non definite) scalar product $<,>$ such that the following conditions hold.
\\
(i) $<[x,y],z> + <y,[x,z]> =0$, $\forall x,y,z\in \Cal G$, where  $[,]$ is the Lie bracket. The pair $(\Cal G,<,> )$ is also said to be an orthogonal Lie algebra.
\\
(ii) $\Cal G= A\oplus B$ possesses two supplementary totally isotropic subalgebras $A$ and $B$ of same dimension $dim (A)=dim(B)$. Vector spaces with this property are sometimes called hyperbolic.

\end{definition}
\begin{definition}
A Left Symmetric Algebra (LSA) structure on a vector space $V$ is a product $V\times V \to V$, $(x,y)\mapsto xy$ whose associator $a(x,y,z):= (xy)z-x(yz)$ is left symmetric i.e $a(x,y,z) = a(y,x,z)$, $\forall x,y,z\in V$. 

 A symplectic Lie group ($G,\omega^+$), in the sense of  Lichnerowicz-Medina {\bf \cite{Li-Me}}, is a
Lie group $G$ together with a left invariant symplectic form $\omega^+$. If
$\Cal G:= Lie(G)$ and  $\omega^+_\epsilon =:\omega$, then $(\Cal G,\omega)$ is
called a symplectic Lie algebra. \end{definition}
An LSA gives rise to a Lie bracket $[,]$ on $V$ by $[x,y]:=xy-yx$.  The Lie backet $[,]$ is also said to be underlying the LSA structure.
Any connected Lie group with Lie algebra $V_-:= (V,[,])$ possesses a left invariant affine structure with associated connection $\nabla_{x^+}y^+:=(xy)^+$. The formula $\nabla_{x^+}y^+_\epsilon =:xy$  gives the converse way to get LSAs on Lie algebras from left invariant affine structures on Lie groups.

In the case of symplectic Lie groups, the connection defined by ({\bf
\ref{PSG}}) induces a Left Symmetric Algebra (LSA) structure $(x,y) \to xy$ in
$\Cal G$, which is compatible with the Lie bracket of  $\Cal G$, ie $xy-yx =
[x,y]$, $\forall x,y\in\Cal G $. The formula {\bf (\ref{PSG})} above, equips $G$
with a left invariant affine structure. 
 In the Lie algebra level, $\omega^+$ being closed reads 

 \begin{equation} \label{co-cyle}
 \omega( [x,y],z)+\omega ([y,z],x) +\omega ([z,x],y)=0
, \hskip 2truemm \forall x,y,z\in \Cal G
 \end{equation}
 On this work, a solution of the
CYBE is always supposed to be skew-symmetric. A left invariant affine structure will always mean a (locally) flat  (and torsion free) left invariant affine structure.

From now on, the scalar product $ <,> $ stands for the duality scalar product (dual pairing) between the considered vector spaces and their duals.



\section{\bf The CYBE and Affine Geometry} 
Here is a simple remark about the Classical Yang Baxter Equation (CYBE), with
several interesting consequences. Denote $Z(\Cal G)$ the center of a Lie
algebra $\Cal G$ .
 \begin{lemme}\label{dualYBC}  Let $(G,\pi)$ be an exact
Poisson  Lie group, i.e $\pi:= r^+ - r^-$, where $r\in \wedge ^2 \Cal G$.
Denote $(\Cal G^*, [,]_*)$ the dual  Lie algebra of $(G,\pi)$. If for every
$\alpha, \beta  $ in $\Cal G^*$ the following assumption holds 
 
 $r (
[\alpha , \beta]_* ) -  [r(\alpha) , r(\beta)] \in Z(\Cal G)$ \hfill $(*)$,
 
 then any Lie group $G^*$ with Lie algebra $(\Cal G^*, [,]_*)$, is endowed
with a left invariant affine connection given by \hskip 1truecm
 $\nabla
_{\alpha^+}\beta^+:= (ad^*_{r(\alpha)}\beta)^+. $
 
 The corresponding LSA
structure on $(\Cal G^*, [,]_*)$  is given by $\alpha, \beta  $ in $\Cal
G^*$\begin{equation} \label{PSGdual} \alpha\beta:= ad^*_{r(\alpha)}\beta.  
\end{equation}  
 \end{lemme}
 
 {\bf Proof.}
  It suffices to check that the map  

 \begin{equation}
\eta : \Cal G^*\to aff(\Cal G^*), \alpha \mapsto (\alpha, ad^*_{r(\alpha)})
\end{equation} is a Lie algebra homomorphism, or equivalently to check the identity 
\begin{equation}\label{Eq-Aff}
([\alpha, \beta]_*,ad^*_{r([\alpha, \beta]_*)} ) = (ad^*_{r(\alpha)}\beta - ad^*_{r(\beta)}\alpha , ad^*_{[r(\alpha), r(\beta)]}). 
\end{equation}
But thanks to $(*)$ and the fact that one has
\begin{equation}\label{dualCrochet} 
[\alpha, \beta]_* = [\alpha, \beta]_r := ad^*_{r(\alpha)}\beta - ad^*_{r(\beta)}\alpha,
\end{equation}
 by definition of the Lie bracket  on the dual of $(G,\pi)$, formula {\bf (\ref{Eq-Aff})} is immediate. \hfill $\square$

\begin{corollaire}\label{corollaire1}Every Lie group is endowed with a left invariant affine structure, if it is dual to a Poisson Lie-group $(G,\pi)$, where $\pi = r^+ - r^-$ is given by a solution $r$ of the CYBE . 
\end{corollaire}

\begin{corollaire}\label{corollaire2} If $G^*$ is a dual Lie group of a Poisson Lie group $(G,\pi)$, where $\pi = r^+ - r^-$ is given by a solution $r$ of the CYBE , then the derived Lie group $[G^*,G^*]$ of $G^*$ is not dense in  $G^*$.
\end{corollaire}
 We can  now state the following result.
\begin{theo}\label{dualYBC-Lag}Let $G$ be a Lie group with (trivialized) cotangent bundle $T^*G= \Cal G^*\rtimes G$, considered as a Lie group (semi-direct product of $G$ and the abelian Lie group $\Cal G^*$, where $G$ acts on $\Cal G^*$ by the coadjoint action).
Then the solutions of the CYBE in $\Cal G := Lie(G)$ are in one-to-one correspondence with the connected Lie subgroups of $T^*G$ equipped with left invariant affine structures whose Lie algebras are lagrangian graphs of a $\Cal G $-valued skew-symmetric linear mapping of $\Cal G^*$. Here  the Lie algebra $t^*\Cal G:= \Cal G^*\rtimes\Cal G$ of $T^*G$, is endowed with its canonical dual pairing $< (\alpha, x), (\beta,y) >:= \alpha (y) + \beta (x)$ for all $x,y\in\Cal G$, $\alpha, \beta \in\Cal G^*$.
\end{theo}
{\bf Proof.} Let $r$  be a solution of the CYBE on $\Cal G$. Let's consider the bi-invariant metric (i.e canonical orthogonal structure) on $T^*G$
  
 $< (\alpha, x)^+, (\beta,y)^+ >:= \alpha (y) + \beta (x)$.

One realizes that the map $\theta : \Cal D (r)\to t^* \Cal G:= Lie (T^* G)$; $\theta (\alpha , x):= (\alpha, r(\alpha) + x)$ is a Manin Lie algebras  isomorphism (where $\Cal D (r)$ is the double Lie algebra of $(G,r^+ - r^-)$), see proposition {\bf \ref{doubleYBC1}}. On the other hand, the Lie  subalgebra $\Cal G^*$ of  $\Cal D(r)$, dual to  $\Cal G$, is equipped with an LSA structure compatible with the Lie bracket of  $(\Cal G^*,[,]_r)=:\Cal G^*(r) $ and defined by $\alpha \beta = ad^*_{r(\alpha)}\beta$. Furthermore,  the image of  $\Cal G^* (r)$  by $\theta$ is  the graph of $r$. Thus, the connected Lie  subgroup of $T^*G$, whose Lie algebra  is $\theta (\Cal G^*(r))= graph (r)$, is indeed endowed with a left invariant affine structure.  

Conversely, let  $K$ be a connected Lie  subgroup of $T^* G$ satisfying the conditions of our statement, in particular $Lie (K) = graph(r)$, where $r:\Cal G^* \to \Cal G$ is skew symmetric. The Lie bracket in $graph (r)$ is of the form 
 $$[\hskip 2truemm (\alpha, r(\alpha)\hskip 1truemm ,\hskip 1truemm  (\beta,r(\beta))\hskip 2truemm ] = ([\alpha, \beta]_r, [r(\alpha), r(\beta)])$$
 where $[\alpha,\beta]_r= [\alpha,\beta]_* = ad^*_{r(\alpha)}\beta - ad^*_{r(\beta)}\alpha$ .

The fact that $graph (r)$ is a subalgebra of $t^* \Cal G$ is equivalent to the following

$ r([\alpha, \beta]_r) = [r(\alpha), r(\beta)]_{\Cal G} $ for all $\alpha, \beta \in \Cal G^*$.

But such an identity (with $r$ skew symmetric,  $graph (r)$ being totally isotropic) means that one has $[r,r] = 0$. In fact, we have $[r,r](\alpha, \beta, \gamma) = <\hskip 1truemm[r(\alpha), r(\beta)] - r([\alpha, \beta]_r)\hskip 1truemm , \hskip 1truemm \gamma \hskip 1truemm>$, for all $\alpha, \beta,\gamma \in\Cal G^*. $  \hfill $\square$

\vskip 5truemm

Let  $G$ be a connected Lie group with Lie algebra $\Cal G$, and  $r\in \wedge^2\Cal G$ a solution of the  CYBE  on $\Cal G$. Now, we would like to study the geodesic completness of the left invariant affine structure induced by $r$ on every Lie group $G^*$ whose Lie algebra is the dual Lie algebra $(\Cal G^*, [,]_r)$ of $(G,r^+-r^-)$. Such a completness occurs if and only if the universal cover $\tilde G^*$ of  $ G^*$ is affinely diffeomorphic to the affine space $\Bbb R^n$, via a developping map. Equivalently, this means that in the LSA $\Cal G^*$,  all the right multiplications have a null trace.  

Consider the left  invariant Poisson tensor $r^+$ on $G$ satisfying $r^+_\epsilon = r$. It's well known that the leaf through $\epsilon$, of the symplectic foliation of $r^+$, is a symplectic Lie group, whose Lie algebra is the image $Im(r)$ of $r:\Cal G^*\to \Cal G $. The natural projection $p:\Cal G^*\to \Cal G^*/ker (r)\simeq Im(r)$ is a morphism of LSAs.

\vskip 2truemm
The symplectic form $\bar \omega$ in $Im(r)$ is given by

 $\bar \omega (r(\alpha), r(\beta)):= <\alpha, r(\beta)>$, $\forall \alpha,\beta\in \Cal G^* $,

 while the LSA structure in  $Im(r)$ is defined as follows

 $\bar \omega (r(\alpha) r(\beta), r(\gamma)):= <ad^*_{r(\alpha)}\beta,r(\gamma)>$ or equivalently $r(\alpha) r(\beta):= r(ad^*_{r(\alpha)}\beta )$.
\vskip 2truemm

 Denoting by $R_\alpha$ (resp. $R^{im}_\alpha$) the right multiplication by $\alpha$ (resp. $r(\alpha)$) in $\Cal G^* $ (resp. in $Im(r)$),
one can show that  $R_\alpha$ and $R^{im}_{r(\alpha)}$ have the same trace. We then deduce the following important result.

\begin{theo}\label{GeoComplAf} ({\bf Geodesic completness of the affine structure})
Let $G^*$ be a Lie group, with Lie algebra $(\Cal G^*,[,]_r)$. The left invariant affine stucture on  $G^*$, induced by $r$, is complete if and only if the leaf $\Cal F_\epsilon$ through $\epsilon$, of the symplectic foliation of $r^+$ is a unimodular Lie group.  So  completness implies the solvability of $G^*$. Here  $r^+$ stands for the left invariant Poisson tensor whose value at $\epsilon$ is $r$.
\end{theo}

Equivalently, we also have\\
\vskip 1truemm
{\bf Theorem \ref{GeoComplAf} bis.}

{\it The left invariant affine stucture on  $G^*$, induced by  $r$, is complete if and only if the quotient $\Cal G^*/Ker (r)$, is a unimodular and solvable Lie algebra. In this case  $G^*$ is solvable, as $Ker (r)$ and $\Cal G^*/Ker (r)$ are. Here $Ker(r)$ is the kernel of the skew symmetric linear map $r:\Cal G^*\to \Cal G$; $\alpha \mapsto r(\alpha,.)= :i_\alpha r$.}\\ 
\vskip 0.5truemm
Remark that $Ker (r)$ is an abelian two-sided ideal of the LSA $\Cal G^*$.


\section{\bf  The double of a  symplectic Lie group is affine and complex}\label{Doubl-alg-sympl}
Denote $\Cal D (\Cal G)$ the double Lie algebra of a Poisson Lie group $(G,\pi)$, where $\pi:= r^+ - r^-$ is given by a solution $r$ of the CYBE on $\Cal G :=Lie(G)$.
Let's point out the following interesting result, whose proof is immediate.
\begin{proposition}\label{doubleYBC1}

The mapping 
 $\Xi: (x,\alpha)\mapsto (x+r(\alpha), \alpha)$  is a Lie algebra isomorphism and an isometry between the Manin algebras $\Cal D( \Cal G)$ and $ \Cal G \ltimes_{ad^*} \Cal G^* $. 
\end{proposition}
 Now, let $(\Cal G , w)$ be a symplectic Lie algebra. Consider the vector space isomorphism 
$$
\begin{matrix}    
 q: \Cal G  \to  \Cal G^* \cr
         x \mapsto q(x) \cr
\end{matrix}
$$
where $<q(x),y> = w(x,y)$ for all  $ x,y$ in $\Cal G$ and 
$<,>$ is the duality scalar product.

Recall that $r:= q^{-1}$ regarded as an element of $\wedge ^2\Cal G$, is a solution of the CYBE on $\Cal G$. Conversely associated to an invertible solution $r$ of the CYBE on $\Cal G$, is a symplectic Lie algebra $(\Cal G,\omega)$ by setting $\omega (x,y):=<r^{-1}(x),y>$.

Now using $q$, let's transport the LSA and the Lie  algebra structures of $\Cal G$ to the vector space $\Cal G^*$, in such a way as $q$ becomes an isomorphism of LSAs. 
One then obtains compatible Lie  algebra and LSA structures on $\Cal G^*$. They are respectively given by
 \begin{equation} \label{YBlin2}
  [\alpha , \beta]_{q}:=q([r(\alpha), r(\beta)])
\end{equation}

 \begin{equation}\label{PSGdual2}
\alpha \beta := q( r(\alpha) r(\beta) )
\end{equation}
for all $\alpha,\beta \in\Cal G^{*}$, where $r(\alpha) r(\beta)$ is the product of $r(\alpha)$ and $r(\beta)$ in  $\Cal G$, induced by {\bf (\ref{PSG})}.

The following can be proved by some direct computations.

 \begin{lemme}\label{doubleYBC2} For every $x, y \in \Cal G$ and $\alpha, \beta \in \Cal G^*$ one has   
$$
\begin{matrix}
(a) \hskip 5truemm  ad^{*}_{x}q(y) = q(xy) =:q(x)q(y) \hskip 3truemm\text{or equivalently } \hskip 3truemm q \circ l_{x} = ad^{*}_{x} \circ q \hskip 2truemm \hfill \cr
(a')\hskip 5truemm r(ad^{*}_{x} \beta) = xr(\beta)= (l_x\circ r)(\beta),\hskip 2truemm r \circ ad^{*}_{x} =(l_x\circ r) \hfill \cr
(b)  \hskip 5truemm  r(ad^{*}_{x} \beta ) =  ad^{*}_{q(x)} r( \beta ) \hfill \cr
(c)  \hskip 5truemm  [ \alpha  ,  \beta ]_{q} =  ad^{*}_{r( \alpha )}  \beta  -  ad^{*}_{r( \beta )} \alpha \hfill \cr
(d)  \hskip 5truemm q(x) \beta = ad^{*}_{x} \beta \hfill \cr
(e)  \hskip 5truemm r( \beta )x = ad^{*}_{ \beta }x \hfill
\end{matrix}
$$
\end{lemme}

\begin{remarques}\label{dual-double}  \end{remarques} 
The formula (d) means that the LSA structure on $\Cal G^*$ defined by the formula {\bf(\ref{PSGdual})} coincides with the one given in {\bf (\ref{PSGdual2})},  which is the image by $q$ of {\bf (\ref{PSG})}.

 Thanks to (c), the Lie bracket $[,]_q$ is the one underlying the LSA structure on $\Cal G^*$ given by  {\bf(\ref{PSGdual})}, and  coincides with the one (namely $[,]_r$ ) defined by {\bf (\ref{dualCrochet})}.
\vskip 3truemm
 Let's consider the coadjoint representations of $\Cal G$ and $(\Cal G^*, [,]_q)$ respectively on $\Cal G^*$ and $\Cal G$. On the vector space $\Cal D:= \Cal  G \oplus \Cal G^* $, one defines the product: $x,y \in \Cal G$ and  $\alpha, \beta \in  \Cal G^* $ 
\begin{equation} \label{double1}  
(x,\alpha)( y,\beta):= (xy + ad^*_\alpha y , \alpha\beta + ad^*_x \beta)
\end{equation}       
 We can now state the main  result of this paragraph, whose proof is given later.
\begin{theo}\label{doubleYBC3} Let $ r\in\wedge^2 \Cal G$ be an invertible solution of the CYBE on a Lie group $G$ with Lie algebra $\Cal G$.
Then every Lie group is endowed with a left invariant affine structure if its  Lie algebra is isomorphic to the double Lie algebra $\Cal D = \Cal G \oplus \Cal G^*$ of the Poisson Lie group $(G, \pi:=r^+ -r^-)$. The associated LSA in $\Cal D$ is given by the formula {\bf (\ref{double1})}. 
\end{theo}

Let $(G,\omega^+)$ be a symplectic Lie group. Consider the left invariant affine connection $\nabla$ in $G$ and the corresponding LSA  structure on  $\Cal G $. 

From theorem {\bf \ref{doubleYBC3}} and proposition {\bf \ref{doubleYBC1} } one deduces the 

\begin{corollaire}\label{doubleYBC4} Let $(G,\omega^+)$ be a symplectic Lie group, with  symplectic Lie algebra $(\Cal G,\omega) $.
 
The map $\Xi$ of proposition {\bf \ref{doubleYBC1}}  gives rise to a left invariant affine (i.e flat and torsion free)  connection $\tilde \nabla$  on the  cotangent bundle $T^*G$ of $G$, in such a way that the natural fiber map (projection) is  a morphism between the affine manifolds $(T^*G,\tilde \nabla )$ and $(G,\nabla)$.
 The LSA structure  induced by $\tilde \nabla$ on $ Lie(T^*G)= \Cal G \ltimes_{ad^*} \Cal G^* $ is given by
 $(x,\alpha)(y, \beta):=(xy, ad^*_x\beta) $, where $xy$ is the product of $x$ and $y$ on $\Cal G $, induced by  {\bf (\ref{PSG})}: $\omega (xy,z):= - \omega (y,[x,z])$, $\forall z\in\Cal G$.
\end{corollaire} 
 Below, we give a method to build algebras with  LSA  structure, that generalizes the theorem {\bf \ref{doubleYBC3}}.

Let  $A$ and $B$ be two LSAs, denote $A_-$ et $B_-$ the respective underlying Lie  algebras.

 \begin{definition} $A$ and $B$  are $(\theta_{1},\theta_{2})$-linked if there exist two Lie algebra representations
 $$ 
\begin{matrix}
 \theta_{1} :A_- & \to & \Cal Gl (B) \hskip 0,5truecm and \hskip 0,5truecm \theta_{2} :B_- & \to & \Cal Gl (A)
\end{matrix}
$$
such that
\begin{equation} \label{comp1}
(\theta_{1}(a)b)b' + b(\theta_{1}(a)b') - \theta_{1}(a)(bb') = \theta_{1} (\theta_{2}(b)a)b' \hskip 0,5truecm and
\end{equation}

\begin{equation} \label{comp2}
(\theta_{2}(b)a)a' + a(\theta_{2}(b)a') - \theta_{2}(b)(aa') = \theta_{2}(\theta_{1}(a)b)a'
\end{equation}
for all elements $a, a'$ in $A$ and $ b, b'$ in $B$. 
\end{definition}
\begin{lemme}\label{doubleASG}

I. Let A and B be two $(\theta_{1},\theta_{2})$-linked LSAs ,
 then the vector space  $V = A \times B$ is endowed with an LSA structure given by the  product of every two elements $(a,b)$ and $(a',b')$ of $A \times B$:  
\begin{equation}\label{Eq1doubleASG}
 (a , b)(a' , b') = (aa' + \theta_{2}(b)a' \hskip 1truemm , \hskip 1truemm bb' + \theta_{1}(a)b') 
\end{equation}
with underlying Lie bracket:
\begin{equation}\label{Eq2doubleASG}
 [(a , b) , (a' , b')] = ( [a , a'] + \theta_{2}(b)a' - \theta_{2}(b')a \hskip 1truemm , \hskip 1truemm [b , b'] + \theta_{1}(a)b' - \theta_{1}(a')b ) 
\end{equation}
II. Conversely every LSA with two  supplementary left ideals, is built with the  above method.
\end{lemme}
Before giving the proof of this lemma, let's remark an interesting example of LSA that can be obtained by such a construction. The LSA corresponding to the  so-called flat left-invariant connections adapted to the automorphism structure of a Lie group, are obtained  in this way, as they arise in {\bf \cite{Me81}} to be a direct sum of two ideals (see {\bf \cite{Me81}} for more details).  

{\bf Proof of lemma \ref{doubleASG}.}
The  bilinearity  is obvious. To prove that {\bf (\ref{Eq1doubleASG})} defines an LSA structure on $V$, it suffices, identifying an element $a$ of $A$ (resp. $b$ of $B$) with the element $(a,0)$ (resp.$(0,b)$) of $ A \times B $; to check that:
$$
\begin{matrix}
(i)   \hskip 5truemm (ab)b' - a(bb') = (ba)b' - b(ab') \cr
(ii)  \hskip 5truemm (ba)a' - b(aa') = (ab)a' - a(ba') \cr
(iii) \hskip 5truemm (bb')a - b(b'a) = (b'b)a - b'(ba) \cr
(iv)  \hskip 5truemm (aa')b - a(a'b) = (a'a)b - a'(ab)
\end{matrix}
$$
Then, one checks that
 
(i) is equivalent to the hypothesis {\bf (\ref{comp1})}

(ii) is equivalent to the hypothesis {\bf (\ref{comp2})}

(iii) expresses the fact $\theta_{2}$ is a representation of the Lie  algebra $B_-$ and

(iv) expresses the fact $\theta_{1}$ is a representation of the Lie algebra $A_-$.
\\
Conversely, let's suppose that $A$ and $B$ are two left (sided) ideals of an LSA $\Cal G = A \oplus B$. 
The maps defined by $\theta_1(a)b:= ab$ and $\theta_2(b)a:= ba$ for all $a\in A$ and all $b\in B$ satisfy the conditions I of this lemma and then enable us to construct an LSA structure on $\Cal G$ which  coincides with the initial one. \hfill  $\square$

\vskip 5truemm
\subsubsection{\bf   Proof of  theorem \ref{doubleYBC3}}
$-$ Proof of the first assertion.

 In the lemma {\bf \ref{doubleASG}}, put $A =  \Cal G$, $B = \Cal G^{*} $  and denote the coadjoint representations of $\Cal G$ and $ \Cal G^{*} $ respectively by $\theta_{1}$ and $ \theta_{2} $. 

For $x\in \Cal G $, $ \alpha $  and  $ \beta $ in $ \Cal G^{*} $ let's consider the following
 
$
 \Theta \hskip 2truemm = \hskip 2truemm ( \theta_{1}(x) \alpha )\beta  +  \alpha ( \theta_{1}(x) \beta )  -  \theta_{1}(x)( \alpha \beta ) \hskip 5truemm 
$
\\
Using the formulas (a)-(e) and the fact that the product defines an LSA structure, one has 
$$
\begin{matrix}
\Theta &=  (ad^{*}_{x} \alpha ) \beta  +  \alpha (ad^{*}_{x} \beta )  -  ad^{*}_{x} ( \alpha \ \beta ) \hfill \cr
&= q( xr( \alpha))q(r( \beta )) + q(r(  \alpha)) q(xr( \beta) ) - q(x(r( \alpha)r( \beta ))) \hfill \cr
&=q((xr( \alpha ) ) r( \beta )  + r( \alpha ) (x r( \beta))  - x(r( \alpha )r( \beta))) \hfill \cr
&= q( ((x r( \alpha )) r( \beta )  - x(r( \alpha )r( \beta )) )  + r(\alpha )(xr( \beta ) )) \hfill \cr
&= q((r( \alpha )x)r( \beta ) - r( \alpha )( x r( \beta )) +  r( \alpha )(xr( \beta ))) \cr
& = q((r( \alpha )x)r( \beta )) \hfill \cr
& = q(ad^{*}_{ \alpha }x)\beta   = ad^{*}_{ad^{*}_{ \alpha } x} \beta \hfill \cr       
& = \theta_{1}( \theta_{2} ( \alpha)x) \beta \hfill 
\end{matrix}
$$
The first hypothesis of lemma {\bf \ref{doubleASG}}   is then satisfied. The second one is exactly checked in the same way.
\\
Then, lemma {\bf \ref{doubleASG}} allows us to make a conclusion to the first statement of this theorem.
\\
$-$The second  statement is directly  deduced from both remarks {\bf \ref{dual-double}} and lemma {\bf \ref{doubleASG}}.\hfill $\square$

\vskip 2truemm
Let $x$ and $\beta$ respectively be elements of the left sided ideals $\Cal G$ and $\Cal G^*$ of $\Cal D (\Cal G)$.  We can show that the trace of the right multiplication by $x$ in $\Cal D (\Cal G)$ is exactly equal to the one of the right multiplication  by $x$ in $\Cal G$, while the trace of the right multiplication by $\alpha$ in $\Cal D (\Cal G)$  is exactly equal to the one of the right multiplication by $r(\alpha)$ in $\Cal G$. Finally we claim

\begin{theo}\label{GeoComplDoubl} ({\bf Geodesic completness of the double}) 

Let $G$ be connected a Lie group with Lie algebra $\Cal G$ and  $r\in \wedge^2\Cal G$ an invertible solution of the CYBE.
 Consider the LSA structure in the  double Lie algebra $\Cal D (\Cal G)$ of $(G,r^+ - r^-)$ as described in the theorem {\bf \ref{doubleYBC3}}.
The following are equivalent: 

1. on every connected Lie group $D (\Cal G)$ whose Lie algebra is $\Cal D (\Cal G)$; the left invariant affine structure corresponding to the LSA, is
geodesically complete. 

2. the Lie group $G$ is unimodular (and solvable). 

\end{theo} 

\vskip 5truemm
Let's denote again by  $\Cal D (\Cal G)$ the double  Lie algebra of a symplectic Lie group $(G,\omega^+)$, and by $q$  as above. One has
 
\begin{proposition}\label{holomorphe} Every Lie group, whose  Lie algebra is the double  Lie algebra $\Cal D (\Cal G)$ of $(G,w)$, is endowed with left invariant complex  structure defined by the endomorphism 
$$
J(x,\alpha) \mapsto (-q^{-1}(\alpha), q(x) )$$ of $\Cal D (\Cal G)$; where  $q(x):= \omega (x,.)$. 
\end{proposition}
To prove the above proposition, one just checks that the Nijenhuis tensor

 $N_J(v_1,v_2)={[Jv_1,Jv_2]-[v_1,v_2]-J[v_1,Jv_2]-J[Jv_1,v_2]}$
\\
 of $J$, vanishes for all $v_1= (x,\alpha), v_2= (y,\beta) \in \Cal D (\Cal G)$.


\vskip 1truecm
\section {\bf The Poisson-Lie structure associated to \\ a symplectic Lie-group ($G,\omega^+$) is polynomial of degree 2}

Let $(G, \omega^+)$ be a symplectic Lie group, $\epsilon$ its unit and $(\Cal G, \omega)$ its symplectic Lie algebra,  where  $\omega =\omega_\epsilon^+ $.

Consider, as in section {\bf \ref{Doubl-alg-sympl}} the vector  spaces (skew-symmetric) isomorphism

$$
\begin{matrix}
q:\Cal G &\to \Cal G^* \hfill \\
x&\mapsto q(x):= i_x \omega = \omega ( x , . ) \hfill
\end{matrix}
$$
Due to the relation $\partial\omega =0$, ie $\omega^+$ is closed, $q$ is actually a 1-cocycle for the coadjoint representation of $\Cal G$. Again, denote $r = q^{-1}$ the inverse of $q$, considered when needed, as an element of $\wedge^2\Cal G$ as well.
\vskip 2truemm
Let $\tilde G$ be the connected and simply connected Lie group with Lie algebra $\Cal G$.
The invertible 1-cocycle $q$ can be integrated into a 1-cocycle $Q: \tilde G \to \Cal G^*$ for the coadjoint representation of $\tilde G$. Such a  1-cocycle $Q$ is a cover over the open subset $ImQ$ of $\Cal G^*$, and is equivariant relative to the actions of $\tilde G$ on itself by left multiplications, and on $\Cal G^*$ by the affine representation

$$
\begin{matrix}
 \tilde G &\buildrel \rho \over \to Aff(\Cal G^*) = \Cal G^* \rtimes GL(\Cal G^*)\\
\sigma &\mapsto \hskip 3truemm \rho (\sigma ) \hskip 5truemm =(Q(\sigma), Ad^*_\sigma)\hfill
\end{matrix}   
$$
Namely, $Q$ is given by the expression
$$
\begin{matrix}
Q(\tau) &= \displaystyle \sum^{k= + \infty }_{k=1} \frac {(ad^*_a)^{k-1}}{k!}q (a) \hskip 9truecm (a)
\end{matrix}
$$
where  $\tau = exp (a)$ and $a$ is in $\Cal G$.

 The affine representation $\rho$ then has an open orbit $Orbit(0) = ImQ$, here $0$ stands for the origine of the vector space $\Cal G^*$.
Pulling back (using $Q$) the usual affine structure on the open subset $Orbit(0)$ of $\Cal G^* $, one obtains a left invariant affine structure on $\tilde G$.
Actually, $Q$ is a developping map of such a structure on $\tilde G$.

\begin{lemme}\label{revet} There exists a left invariant affine structure on $G$ such that the ``universal covering map'' $\tilde G \buildrel p \over \to G$ is a morphism of affine manifolds. Such an affine structure coincides with the one given by $\omega ^+$ on $G$. Furthermore $p$ is Poisson morphism between the Poisson Lie-groups $(\tilde G, r^+ - r^+)$ and  $(G,r^+ - r^+)$.  
\end{lemme}  
{\bf Proof.}

As the group of deck transformations $\pi_1 (G)$ is a  central and discrete subgroup of $\tilde G$, it  properly and freely acts on  $\tilde G$ by  left translations which are affine transformations.
 Consequently, $G = \tilde G \slash {\pi_1 (G) }$ is equipped with a left invariant affine structure relative to which, $p$ is a morphism of affine manifolds.
\\
Realizing that $q$ is the  orbital map (for the infinitesimal action) at $0$, we deduce the coincidence between such an affine structure on $G$ and the one given by  $\omega^+ $. 
\\
Let's denote $\tilde \pi$ (resp. $\pi$), the tensor on  $\tilde G$ (resp. $G$) given by $r$. Let $\tilde \epsilon$ be the unit of  $\tilde G$ and $\tilde L_\sigma$ (resp. $ \tilde  R_\sigma$) be the left (resp. right) translation  in $\tilde G$ by $\sigma \in \tilde G$.
Then, for every $\sigma \in \tilde G$, one has
\\
$T_\sigma p \hskip 1truemm\tilde \pi_\sigma = T_\sigma p (T_{\tilde\epsilon}\tilde L_\sigma\hskip 1truemm r) - T_\sigma p (T_{\tilde\epsilon}\tilde  R_\sigma \hskip 1truemm r) = T_{\tilde\epsilon}( p \circ\tilde L_\sigma)\hskip 1truemm r - T_{\tilde\epsilon } (p \circ\tilde R_\sigma)\hskip 1truemm r = (T_\epsilon  L_{p(\sigma)})\hskip 1truemm r - (T_\epsilon  R_{p(\sigma)})\hskip 1truemmr =: \pi_{p(\sigma)} $.

Then, $p$ is a Poisson morphism between $(\tilde G, \tilde \pi)$ and $(G,\pi  )$.\hfill $\square$
\begin{remarques}\label{revet-aff}

1. As the covering map  $p$ is a morphism  of   both affine and Poisson manifolds,  to proof that $\pi$ is polynomial, it suffices to prove that $\tilde\pi$ is polynomial.

2. It's convenient to point out that for every $\gamma \in \pi_1 (G)$,  the following diagram commutates
$$
\begin{matrix}
\tilde G &\buildrel Q \over \to &\Cal G ^* \\
 \gamma \downarrow && \downarrow h(\gamma) \hfill \\
  \tilde G &\buildrel Q \over \to &\Cal G^* 
 \end{matrix}
$$

Here $h: \pi_1 (G) \to Aff(\Cal G^*)$ is the Lie groups homomorphism, the so called holonomy representation of the affine  manifold $G$.

3. One easily verifies that $Q$ is a moment map for the action of $\tilde G$ on itself, via left multiplications. 
\end{remarques}
 
Here is our main result of this section. It states that for every invertible solution $r$ of the CYBE, the corresponding Poisson Lie structure $\pi:=r^+-r^-$ is polynomial, of degree smaller or equal to $2$.

\begin{theo}\label{local-aff1} Let $(G,\omega ^+)$ be a connected symplectic Lie group.  Denote $r$ the associated solution of CYBE and  $r^+$ (resp.$r^-$) the corresponding left (resp. right) invariant Poisson tensor.

Then, the Poisson Lie tensor $\pi = r^+ - r^-$ is polynomial, of degree smaller or equal to $2$, relative  to the left invariant affine structure on  G defined by $\omega ^+$. That is, at the neighborhood of every point of $G$, there exists local  affine coordinates $(x_1, ..., x_{2n})$ satisfying

$\{x_i, x_j\}_\pi= \displaystyle \sum_s b_sx_s + \displaystyle \sum_{l,k} b_{lk}x_lx_k$ where the $b_s, b_{lk}$'s are in $\Bbb R$. 
\end{theo}
\subsection{\bf Proof of theorem \ref{local-aff1}}

To prove that $\pi$ is polynomial, thanks to remarks {\bf \ref{revet-aff}}, it suffices to prove that $\tilde \pi$ is polynomial.
 We can then suppose,  without any lost of generalities, that $G$ is connected and  simply connected. 

We will then prove that $\pi$ can be projected by $Q$ to  a Poisson tensor on $ImQ$ which is polynomial, of degree smaller or equal to $2$. 

The following very useful two lemmas will be proved later.

\begin{lemme}\label{local-aff2}

 Let $\sigma$ be in $G$ and $T_\sigma Q$ be the differential of $Q$ at $\sigma$, 
the tensor $T_\sigma Q$ $\pi _\sigma$ only  depends on $Q(\sigma)$.
\end{lemme} This means that if two elements $\sigma$ and $\sigma '$ of $G$ have the same image $Q(\sigma) =Q(\sigma') $ by $Q$,  then one has 
 $T_\sigma Q$ $\pi _\sigma = T_{\sigma '} Q$ $\pi _{\sigma '}$.

The lemma {\bf \ref{local-aff2}} solves the main problem and then proves that the  tensor $\pi$ can be projected, by $Q$, to  a Poisson tensor $\Lambda = Q_*\pi$ on $ImQ$. 
\begin{lemme}\label{local-aff3}
Let $f$ and $g$ be two elements of $ (\Cal G^*)^* = \Cal G$, denote $\{f,g\}_\Lambda$ their Poisson brackets relative to $\Lambda$. Then one has
$$
\begin{matrix}
\{f,g\}_\Lambda (Q(\sigma))& =- < Q(\sigma)\hskip 1.5truemm,\hskip 1.5truemm [f,g]> -  r (ad^*_fQ(\sigma)\hskip 1.5truemm,\hskip 1.5truemm ad^*_g Q(\sigma)) \hfill \\
&= - < Q(\sigma)\hskip 1.5truemm,\hskip 1.5truemm [f,g]> - <(ad_f \otimes ad_g ) r\hskip 1.5truemm,\hskip 1.5truemm  Q(\sigma)\otimes Q(\sigma) > \hfill
\end{matrix}
$$
\end{lemme} Then, the bracket of the restriction to $ImQ$ of two linear functions on $\Cal G^*$, is a polynomial function of degree $2$, on $ImQ$. 
Being understood that $[f,g]$ is a linear form on $\Cal G ^*$, the function 
$$
\begin{matrix}
u: ImQ \to \Bbb R \hfill \\ 
\hskip 6truemm Q(\sigma) \mapsto u(Q(\sigma)):=   r (ad^*_fQ(\sigma), ad^*_g Q(\sigma)) = <(ad_f \otimes ad_g ) r,  Q(\sigma)\otimes Q(\sigma) >
\end{matrix}
$$
is polynomial (with respect to the coordinates of $Q(\sigma) )$ homogeneous of  degree $2$, and the element
$(ad_f \otimes ad_g ) r $ of $\otimes^2 \Cal G$ is a bilinear form on $ \Cal G^*$.


This ends the proof of the theorem. \hfill $\square$
\vskip 5truemm

\subsection{ \bf Proof of lemmas \ref{local-aff2} et \ref{local-aff3}}
\vskip 3mm
{\bf a. Proof of lemma  \ref{local-aff2} }

Let $f$ and $g$ be in  $\Cal C^\infty (ImQ)$; one has
$$
\begin{matrix}
T_\sigma Q\hskip 0.5truemm \pi_\sigma ( df(Q(\sigma))\hskip 1.5truemm,\hskip 1.5truemm dg(Q(\sigma)))& = \pi _\sigma (df(Q(\sigma))\circ T_\sigma Q \hskip 1.5truemm,\hskip 1.5truemm dg (Q(\sigma))\circ T_\sigma Q) \hfill \\
&= < \pi^\sharp _\sigma (df(Q(\sigma))\circ T_\sigma Q)\hskip 1.5truemm ,\hskip 1.5truemm dg (Q(\sigma))\circ T_\sigma Q > \hfill\\
& = <(T_\sigma Q \circ \pi^\sharp _\sigma \circ (T_\sigma Q)^t)\hskip 1truemm (df(Q(\sigma))\hskip 1.5truemm ,\hskip 1.5truemm dg (Q(\sigma))> \hskip 26truemm (a)
\end{matrix}
$$  
 Let's then prove that the expression $T_\sigma Q \circ \pi^\sharp _\sigma \circ (T_\sigma Q)^t$ only depends on  $Q(\sigma)$. 
To do this, we will need the following formulas.
Denote by $\Omega$ the differential 2-form  $\Omega = \omega ^+ - \omega ^-$ and $\Omega ^b$ the corresponding morphism of vector bundles: 
$$
\begin{matrix}
 T_\sigma G \buildrel \Omega^b _\sigma \over \to T^*_\sigma G \hfill\\
 X_\sigma \mapsto \Omega_\sigma^b (X_\sigma) \hfill \\
\end{matrix}
$$ where $<\Omega ^b _\sigma (X_\sigma) , Y_\sigma > := \Omega_\sigma ( X_\sigma, Y_\sigma)$  for all
$X_\sigma, Y_\sigma \in  T_\sigma G$. 

Then for all $\sigma \in G$, the equality 
$$
\begin{matrix}
^t(T_\epsilon L_\sigma ) \circ \Omega ^b _\sigma \circ T_\epsilon L_\sigma = q - Ad^*(\sigma ^{-1}) \circ q \circ Ad(\sigma)  \hskip 7truecm (b)
\end{matrix}
$$ is due to the following, satisfied by every
 $x$ and $y$ in $\Cal G$:   
$$
\begin{matrix}
 <(q - Ad^*_{\sigma ^{-1}} \circ q \circ Ad_\sigma) (x) \hskip 1truemm,\hskip 1truemm y >& = w(x,y) - w(Ad_\sigma x\hskip 1truemm,\hskip 1truemm Ad_\sigma y) \hfill \\
&= \omega^+ _\sigma (x^+_\sigma , y^+_\sigma) - \omega^-_\sigma (x^+_\sigma , y^+_\sigma)\hfill \\ & = \Omega_\sigma (x^+_\sigma , y^+_\sigma) = <\Omega_\sigma ^b (x^+_\sigma ), y^+_\sigma> \hfill \\
& = <^t(T_\epsilon L_\sigma)(\Omega ^b _\sigma \circ T_\epsilon L_\sigma) (x)\hskip 1truemm,\hskip 1truemm y > .
\end{matrix}
$$ 
For $\sigma \in G$  the map $q^\sigma =   Ad^*_{\sigma ^{-1}} \circ q \circ Ad_\sigma$ is a 1-cocycle of $\Cal G$. This is due to the fact that  $Ad_\sigma$ is a Lie algebra homomorphism and $q$ is a 1-cocycle of $\Cal G$.

 Let's integrate $q^\sigma$ into a 1-cocyle  $Q^\sigma: G\to \Cal G^*$, for the coadjoint action of $G$.
Set $\tau = exp (a)$,  where $a$ is in $\Cal G$. Then
$$
\begin{matrix}
Q^\sigma (\tau) &= \displaystyle \sum^{k=+\infty}_{k=1} \frac {(ad^*_a)^{k-1}}{k!} (Ad^*_{\sigma ^{-1}} \circ q \circ Ad_\sigma)(a)\hfill \\
&= Ad^*_{\sigma^{-1}} \displaystyle \sum^{k=+\infty}_{k=1}\frac{(ad^*_{(Ad_\sigma a)})^{k-1}}{k!} (q (Ad_\sigma a)) \hfill \hskip 8truecm (c)
\end{matrix}
$$ 
as $ad^*_a \circ Ad^*_{\sigma{-1}}$ is equal to $Ad^*_{\sigma ^{-1}}\circ ad^*_{Ad_\sigma a}$, for every  $a$ in $\Cal G$ and $\sigma$ in $G$. The last assertion is due to the following
 $$
\begin{matrix}
< ad^*_a \circ Ad^*_{\sigma{-1}}\alpha \hskip 1truemm,\hskip 1truemm x > &= - <\alpha \hskip 1truemm,\hskip 1truemm Ad_\sigma [a,x] > =- <\alpha \hskip 1truemm ,\hskip 1truemm [Ad_\sigma a, Ad_\sigma x ]> \hfill \\
&= <Ad^*_{\sigma ^{-1}} \circ ad^*_{Ad_\sigma a}\alpha \hskip 1truemm,\hskip 1truemm x>\hfill
\end{matrix}
$$ for all $x$ in $\Cal G$ and $\alpha$ in $\Cal G^*$. 
 $$
\begin{matrix}
Q^\sigma (\tau)&= Ad^*_{\sigma^{-1}} (Q(exp(Ad_\sigma a))) = Ad^*_{\sigma^{-1}} Q(\sigma \tau \sigma^{-1})\hfill \\
& = Ad^*_{\sigma^{-1}}( Q(\sigma)+ Ad^*_\sigma Q(\tau \sigma^{-1}))\hfill \\
& = Ad^*_{\sigma^{-1}}Q(\sigma)+ Q(\tau) + Ad^*_\tau Q(\sigma^{-1}))\hfill \\
& = - Q(\sigma^{-1})+ Q(\tau) + Ad^*_\tau Q(\sigma^{-1}))\hfill \\
\end{matrix}
$$ 
Because, one has
\\
 $exp(Ad_\sigma a) = \sigma exp(a)\sigma^{-1}$ = $\sigma \tau \sigma^{-1}$ and $Q(\sigma ^{-1}) = - Ad^*_{\sigma ^{-1}} Q(\sigma)$.
\\
 We then obtain the equality \\
 $ ( Q -Q^\sigma )(\tau) =   Q(\sigma^{-1}) - Ad^*_\tau Q(\sigma^{-1})) = - d(Q(\sigma^{-1}))(\tau).$ \hfill  (d)
\\
Here $d$ is the (Chevalley-Eilenberg) differential associated to the coadjoint representation of $G$ by
 $d\alpha (\sigma): = Ad^*_\sigma \alpha - \alpha$, if $\alpha \in \Cal G^*$  and $ \sigma \in G$.

 In other words, $Q$ and $Q^\sigma$ are cohomologous and their difference is the 1-coboundary $Q - Q^\sigma = - d(Q(\sigma^{-1}))$. Differentiating at $\epsilon$ the relation (d), one establishes that $q$ and $q^\sigma$ are also
cohomologous, and 
\\
$q - q^\sigma = -\delta (Q(\sigma^{-1}))$ \hfill (e)
\\
Here, $\delta$ is the differential associated to (infinitesimal) coadjoint action of $\Cal G$: if $x\in \Cal G$ and $\alpha \in \Cal G^*$ we have
$\delta \alpha (x) := ad^*_x \alpha$.
\\
The relation (b) above, then reads 
\\
$ \Omega ^b = - (T_{\sigma}L_{\sigma {-1}})^t\circ \delta (Q(\sigma^{-1}))\circ T_{\sigma}L_{\sigma {-1}}$ \hfill (b')
\\
It is easy to check that one has 

($r^+)^\sharp = T_\epsilon L_\sigma \circ r \circ (T_\epsilon L_\sigma)^t$
;  $(r^-)^\sharp=T_\epsilon R_\sigma \circ r \circ (T_\epsilon R_\sigma)^t$ \hfill (f)
\vskip 1truemm
$q_\sigma ^+:= (\omega ^+)^b_\sigma = (T_\sigma L_{\sigma^{-1}})^t \circ q \circ (T_\sigma L_\sigma^{-1})$;

 $ q_\sigma ^-:= (\omega ^-)^b_\sigma = (T_\sigma R_\sigma^{-1})^t \circ q \circ (T_\sigma R_\sigma^{-1})$ \hfill (g)
\\
We also have $T_\sigma Q = Ad^*_\sigma \circ q \circ T_\sigma L_{\sigma {-1}}, $ \hfill (h)
\\
 as the following equalities are true, for every
 $x$ in $\Cal G$
$$
\begin{matrix}
 T_\sigma Q (x^+_\sigma ) &= \frac{d}{dt}|_{t=0} Q(\sigma exptx ) = \frac{d}{dt}|_{t=0} (Q(\sigma ) + Ad^*_\sigma Q(exptx)) \hfill \\ 
&= Ad^*_\sigma q(x) = Ad^*_\sigma \circ q \circ T_\sigma L_{\sigma ^{-1}}(x^+_\sigma) . \hfill \\
\end{matrix}
$$
In such  a way that we have 
$$
\begin{matrix}
^t(T_\sigma Q) &= (T_\sigma L_{\sigma ^{-1}})^t \circ q^t  \circ Ad_{\sigma {-1}}=- (T_\sigma R_{\sigma ^{-1}})^t \circ Ad^*_\sigma \circ q  \circ Ad_{\sigma {-1}}\hfill \\
&=- (T_\sigma R_{\sigma ^{-1}})^t \circ q^{\sigma ^{-1}} = - (T_\sigma R_{\sigma ^{-1}})^t \circ (q + \delta (Q(\sigma)))\hskip 6truecm (i)
\end{matrix}
$$
From equality $\Omega^b = q^+ - q^-$ and, from the fact that the morphism $(\tilde r^+)^\sharp \circ q^+$ (resp. $ q^-\circ (\tilde r^- )^\sharp$ ) is nothing but the identity map of $T_\sigma G$ (resp. $T^*_\sigma G$), at each point $\sigma$ of $G$, we get 


$\pi^\sharp = - (\tilde r ^+)^\sharp \circ \Omega^b \circ (\tilde r^-)^\sharp$ \hfill (j) 
 \\
 Noticing the relation $ Q(\sigma) = -Ad^*_\sigma Q(\sigma^{-1})$ first we have for each $x\in\Cal G$ 
$$
\begin{matrix}
\delta (Q(\sigma))(x):= ad^*_xQ(\sigma)
 &= - (ad^*_x \circ Ad^*_\sigma )Q(\sigma ^{-1})
\cr
 &= - (Ad^*_\sigma \circ ad^*_{(Ad_{\sigma^{-1}} x)} )Q(\sigma^{-1}).
\end{matrix}
$$
Second, we obtain 
$$
\begin{matrix}
(Ad^*_\sigma \circ \delta (Q(\sigma^{-1})) \circ Ad_{\sigma ^{-1}})(x) &= Ad^*_\sigma (\delta (Q(\sigma^{-1})) (Ad_{\sigma ^{-1}}x )) 
\cr
 &=Ad^*_\sigma \circ  ad^*_{(Ad_{\sigma ^{-1}x)}}(Q(\sigma^{-1})) .
\end{matrix}
$$
This proves the following 
$$
\begin{matrix}
\delta (Q(\sigma)) = - Ad^*_\sigma \circ \delta (Q (\sigma ^{-1})) \circ Ad_{\sigma ^{-1}} \hskip 9.9truecm (k)
\end{matrix}
$$
Now, let's prove, using formulas above, that

  $T_\sigma Q \circ \pi^\sharp _\sigma \circ (T_\sigma Q)^t$
\\
 only depends on $Q(\sigma)$.
Formulas (h), (i) and (j) enable to write 
$$
\begin{matrix}
 T_\sigma Q \circ \pi^\sharp _\sigma \circ (T_\sigma Q)^t &= Ad^*_\sigma \circ q \circ T_\sigma L_{\sigma {-1}}\circ (\tilde r^+)^\sharp \circ \Omega^b \circ (\tilde r^-)^\sharp \hfill\\
&\hskip 5truemm \circ (T_\sigma R_{\sigma ^{-1}})^t \circ (q + \delta (Q(\sigma))) \hfill 
\end{matrix} 
$$
 Apply now  (b') and (f)  to get 
$$
\begin{matrix}
 T_\sigma Q \circ \pi^\sharp _\sigma \circ (T_\sigma Q)^t &=- Ad^*_\sigma \circ q \circ T_\sigma L_{\sigma {-1}}\circ (T_\epsilon L_\sigma \circ r \circ (T_\epsilon L_\sigma)^t ) \hfill \\
&\hskip 5truemm \circ ^t(T_{\sigma}L_{\sigma {-1}})\circ \delta (Q(\sigma^{-1}))\circ T_{\sigma}L_{\sigma {-1}}    \hfill \\
&\hskip 5truemm\circ ( T_\epsilon R_\sigma \circ r \circ (T_\epsilon R_\sigma)^t ) \circ (T_\sigma R_{\sigma ^{-1}})^t \circ (q + \delta (Q(\sigma))) \hfill\\
\end{matrix}
$$
Simplifying this, one obtains
$$
\begin{matrix}
 T_\sigma Q \circ \pi^\sharp _\sigma \circ (T_\sigma Q)^t =- Ad^*_\sigma \circ
 \delta (Q(\sigma^{-1}))\circ Ad_{\sigma {-1}} \circ r\circ (q + \delta (Q(\sigma))) 
\end{matrix}
$$

 Last, the equality (k) allows us to conclude that the wanted expression (a) reads 

$
 T_\sigma Q \circ \pi^\sharp _\sigma \circ (T_\sigma Q)^t = \delta (Q(\sigma))\circ r\circ (q + \delta (Q(\sigma)))$ \hfill (l)
\\
and then only depends on $Q(\sigma)$. This proves lemma {\bf \ref{local-aff2} }. \hfill $\square$
\vskip 3truemm
{\bf b. Proof of lemma \ref{local-aff3} }

From lemma {\bf \ref{local-aff2} }, the tensor $\pi$ can be projected by $Q$, into a Poisson tensor $\Lambda =Q_*\pi$ on $ImQ$, whose associated vector bundles morphism is
$\Lambda^\sharp_{Q(\sigma)} = T_\sigma Q \circ \pi^\sharp _\sigma \circ (T_\sigma Q)^t$ on the tangent fiber over any point $Q(\sigma) $ of $Im Q$.
Let $f$ and $g$ be in $\Cal C^\infty (ImQ)$, if $\{f,g\}_\Lambda$ is the Poisson bracket of $f$ and $g$; using expression (l) and the definition of $\delta$, we establish the following.
$$
\begin{matrix}
\{f,g\}_\Lambda (Q(\sigma))& = <(\Lambda)^\sharp (df(Q(\sigma))), dg(Q(\sigma)) >\hfill \\
 &= < \delta (Q(\sigma)) \hskip 1truemm (r \hskip 1truemm (q + \delta (Q(\sigma)) \hskip 1truemm)\hskip 1truemm (df\hskip 1truemm(Q(\sigma))) \hskip 1truemm , \hskip 1truemm dg(Q(\sigma)) > \hfill \\
 &=  < \delta (Q(\sigma )) \hskip 1truemm (df(Q(\sigma))), dg(Q(\sigma)) > \hfill \\
& \hskip 4truemm +  < \delta (Q(\sigma)) \hskip 1truemm (r\circ \delta (Q(\sigma))) (df(Q(\sigma))), dg(Q(\sigma)) > \hfill \\
  &=- <Q(\sigma )\hskip 1truemm , \hskip 1truemm [x, y] > 
+  < \delta (Q(\sigma)) \hskip 1truemm (r (ad^*_x Q(\sigma)), y> \hfill \\ 
&=- <Q(\sigma )\hskip 1truemm , \hskip 1truemm [x, y] > 
- < \delta (Q(\sigma)) (y)\hskip 1truemm,\hskip 1truemm r (ad^*_x Q(\sigma))> \hfill \\
&=- <Q(\sigma )\hskip 1truemm , \hskip 1truemm [x, y] > 
- r (ad^*_x Q(\sigma\hskip 1truemm),\hskip 1truemm ad^*_y Q(\sigma ) \hfill \\
&=- <Q(\sigma )\hskip 1truemm , \hskip 1truemm [x, y] > - r (ad^*_f Q(\sigma\hskip 1truemm),\hskip 1truemm ad^*_g Q(\sigma ) \hfill \\
&=-<Q(\sigma )\hskip 1truemm , \hskip 1truemm [x,y ]> - < (ad_x \otimes ad_y) r , Q(\sigma ) \otimes Q(\sigma) >. \hfill 
\end{matrix}
$$  
 Here  $df(Q(\sigma ))=:x$  and $dg(Q(\sigma))=:y$ are seen as elements of the Lie algebra $\Cal G$.

If $f$  and  $g$  are the restrictions  to $Im Q$ of some linear functions (also denoted $f$ and $g$ respectively) on $\Cal G^*$, we have $df(Q(\sigma) = f =: x$, $dg(Q(\sigma) =g = :y$ and the two last equalities above then read
  $$
\begin{matrix}
\{f,g\}_{Q_*(\pi)}(Q(\sigma))& =- <Q(\sigma )\hskip 1truemm , \hskip 1truemm [f, g]> - r (ad^*_f Q(\sigma\hskip 1truemm),\hskip 1truemm ad^*_g Q(\sigma ) \hfill \\
&=-<Q(\sigma )\hskip 1truemm , \hskip 1truemm [f, g]> - < (ad_f \otimes ad_g) r , Q(\sigma ) \otimes Q(\sigma) >. \hfill 
\end{matrix}
$$  
This proves lemma {\bf \ref{local-aff3}}. \hfill $\square$

To end this section, let's point out a consequence of the above theorem {\bf \ref{local-aff1}}. We will use the same notations as in theorem {\bf \ref{local-aff1}}. 
Let $F_\sigma$  be the symplectic leaf,  through $\sigma \in G$, of the Poisson Lie group $(G,\pi:= r^+ - r^-)$. The
vector space $T_\sigma F_\sigma$ tangent to $F_\sigma$, at the point $\sigma$,  is given by 
$T_\sigma F_\sigma = T_\epsilon R_\sigma\{ \eta (\sigma)(Ad^*_\sigma \alpha , \hskip 1truemm .), \alpha \in \Cal G^*\}$.  Here we have $\eta (\sigma): = Ad_\sigma \tilde r - \tilde r$ and $R_\sigma$ is the right multiplication by $\sigma$ in $G$.

 Notice the following relation $\eta (\sigma)(Ad^*_\sigma\alpha , \hskip 1truemm .) = Ad_\sigma (r(Ad^*_\sigma\alpha))  $. Due to formula {\bf (h)}, one has $ T_\sigma Q (T_\sigma F_\sigma) = \{ ad^*_{(r\circ Ad^*_\sigma)\alpha}Q(\sigma),  \hskip 2truemm \alpha \in \Cal G^*\} $.

 As  $r\circ Ad^*_\sigma: \Cal G^* \to \Cal G$ is a vector spaces isomorphism, $T_\sigma Q$ bijectively maps $T_\sigma F_\sigma$ onto the orbit of $Q(\sigma)$ for the $\Cal G$-coadjoint action. 

Namely, we have

\begin{corollaire} \label{feuille-coadj} Let $(G,\omega ^+)$ be a connected and simply connected symplectic Lie group. Put $q(x):=i_x\omega$, $x\in \Cal G:=Lie(G)$ and $r = q^{-1}$, define $Q$ as above.

The leaves of (the symplectic foliation of) the Poisson Lie group $(G,\pi = r^+ - r^-)$ are coverings over the traces on $ImQ$ of the  coadjoint orbits of $G$.  In particular, if $\omega ^+$ is exact, those symplectic   leaves are coverings over coadjoint orbits of $G$.

There is a bijection between the  Lie subgroup  $G_\pi= \{\sigma \in G, \pi_\sigma = 0\}$ of $G$, and the intersection of $ImQ$ and the  vector subspace $\{ \alpha \in \Cal G^*$ s.t. $\alpha$ is invariant relative to the coadjoint representation$\}$ of $G$.    

\end{corollaire}


{\footnotesize
Departement of mathematics. National University of Ireland, Maynooth. Co. Kildare, Ireland.
\\
UMR 5030 du CNRS. D\'epartement de Math\'ematiques.\\ 
Universit\'e Montpellier 2, 34095 Montpellier cedex 5.
\\ 
e.mail adresses: adiatta@maths.may.ie; or  andrediatta@hotmail.com;\\ medina@darboux.math.univ-montp2.fr }
\end{document}